\documentclass[10pt, oneside, draft]{amsart}

        \usepackage {amssymb,amsthm,amsxtra,amscd,amsmath,amsfonts}
        
        \newcommand{\K}{\ensuremath{\Bbbk}}
        \newcommand{\Q}{\ensuremath{\mathbb{Q}}}
        \newcommand{\Z}{\ensuremath{\mathbb{Z}}}
        \newcommand{\CG}{\ensuremath{\mathcal{G}}}
        \newcommand{\CA}{\ensuremath{A}}

        \newcommand{\slin}{sl_2}
        \newcommand{\D}{D_q^{\Lambda}(R)}
        \newcommand{\grHom}{\textup{grHom}}
        \newcommand{\grhom}{\grHom_\K(R, R)}
        \newcommand{\tensor}{\bigotimes}
        \newcommand{\dbox}[1]{[\![#1]\!]}
        \newcommand{\mb}[1]{\mathbf{#1}}
        \newcommand{\mba}{{(1,1)}}
        \newcommand{\mbb}{{(-1,1)}}

        \theoremstyle{plain}
                \newtheorem{theorem}{Theorem}[subsection]
                \newtheorem{lemma}{Lemma}[subsection]
                \newtheorem{corollary}{Corollary}[subsection]
                \newtheorem{proposition}{Proposition}[subsection]
                \newtheorem{remark}{Remark}[subsection]
                \newtheorem*{definition}{Definition}
        \newcommand{\term}[1]{\textit{#1}}

\begin{document}
\setcounter{section}{-1}


        \title[$q$-differential operators]{Quantum differential operators on 
                        the quantum plane}
        \date{\today}
        \author[Iyer]{Uma N. Iyer} 
        \email{uiyer@mri.ernet.in}
        \author[McCune]{Timothy C. McCune} 
        \email{tim@mri.ernet.in}
        \address{The Mehta Research Institute, Chhatnag Road, Jhusi\\
                  Allahabad, U.P. 211 019\\
                  India}
        \begin{abstract}
                The universal enveloping algebra $U(\CG)$ of 
                a Lie algebra $\CG$ acts on its representation
                ring $R$ through $D(R)$, the ring of differential operators
                on $R$.  A quantised universal enveloping algebra (or 
                \textit{quantum group}) is a deformation of a universal
                enveloping algebra and acts not through the 
                differential operators of its representation ring but
                through the quantised differential operators of its 
                representation ring.  We present this situation for the 
                quantum group of $sl_2$.
        \end{abstract}
\maketitle              


\section{Introduction}
   Let $q$ be a transcendental element over $\Q$, and let $\K$ be a field 
   extension of $\Q(q)$ containing $\sqrt{q}$. Let $U_q$ denote the
   quantum group corresponding to the Lie algebra $sl_2(\K )$.
   Let
   $
        R = \K \langle x,y \rangle \diagup xy - qyx.
   $
   We will call $R$ the \term{coordinate ring of the quantum plane} 
   or sometimes just the \term{quantum plane}.
   This ring $R$ is a \textit{representation ring} of $U_q$;  that is, 
   every \textit{type-1}, irreducible, finite dimensional representation
   of $U_q$ appears in $R$ exactly once. Hence,
   $U_q$ acts on the quantum plane. 
   This action is through the quantum- (or $q$-) differential operators 
   (Section~3.3 of \cite{LR1}).  The weight space of $U_q$ is $\Z$,
   the group of integers.  Thus $R$ is $\Gamma = \Z \times \Z$-graded as
   $\deg(x) = (1,1)$ and $\deg(y) = (-1,1)$.  This corresponds
   to the fact that $x$ (resp. $y$)
   can be seen as the highest (resp. lowest) weight vector of weight
   1 (resp. -1) of the unique, type-1, simple, 2-dimensional module.

   The ring $R$ is also graded by the subgroup $\Lambda $ generated by
   $(1,1)$ and $(-1,1)$.  In
   this paper we compute the ring (denoted by $D^{\Lambda}_q(R)$)
   of $q$-differential operators of
   $R$ viewing $R$ as $\Lambda$-graded, and study its properties. 
   Note that the ring $D_q^{\Lambda}(R)$ differs from the ring
   $D_q^{\Gamma}(R)$ of $q$-differential operators on $\Gamma$-graded $R$.
   We shall address this in Section~\ref{S:D0}.

   In Section~\ref{S:D0}, we give the requisite preliminaries.
   Section~\ref{S:D1} deals with the description of first order 
   $q$-differential operators on $\Gamma$-graded $R$.
   In Section~\ref{S:Dlambda} we find the generators of $D^{\Lambda}_q(R)$ 
   explicitly.
   In Section~\ref{S:pro} 
   we describe some basic properties of $D^{\Lambda}_q(R)$
   as a ring.  In particular, we show that $D^{\Lambda}_q(R)$ is isomorphic
   to $D_q(\K [t])\bigotimes _{\K}D_q(\K [t])$ as rings (the ring
   $D_q(\K [t])$ has been studied in \cite{IM}) which is simple by
   Proposition \ref{P:tensor}.  We also show that $D^{\Lambda}_q(R)$ 
   is a domain.
   In Section~\ref{S:mot} we explain the relevance of $U_q$ with 
   $D_q^{\Lambda}(R)$.
   As an application, we show that 
   any irreducible, finite dimensional 
   $U(sl_2(\K))$-module can be given a type-1, simple $U_q$-module structure 
   where $U(sl_2(\K))$ denotes the enveloping
   algebra of $sl_2(\K)$ (Corollary~\ref{C:UtoU_q}). 

   We would like to acknowledge the formulae derived in the Mathematical
   Physics literature (see \cite{DP}, \cite{D}), 
   which have been rederived here.
\section{Some definitions}\label{S:D0}
  In this section we recall some basic definitions from Section~3 of 
  \cite{LR1}.  Let $G$ be an abelian
  group, and $S$ be a $G$-graded $\K$-algebra.  We fix a bicharacter 
  $\beta :G \times G \to \K^*$. 
  The ring of $G$-graded and $\K$-linear
  homomorphisms, $\grHom_{\K}(S,S)$ is an $S$-bimodule given by
  $(\varphi \cdot r)(s) = \varphi (rs)$ and 
  $(r \cdot \varphi)(s) = r\varphi(s)$.
  For homogeneous $r$ of degree $a$,
  we let $[\varphi ,r]_g = \varphi \cdot r - \beta (g, a) r\cdot \varphi$.
  For $g = 0$ the identity element, we write $[\varphi ,r]_0$ simply as
  $[\varphi ,r]$.  This notation can be immediately 
  generalised for a homogenous $\psi$ of degree $a$ 
  (and hence for any graded $\psi$) 
  as $[\varphi ,\psi ]_g = \varphi \psi - \beta (g,a)
   \psi \varphi$.
   Again we let $[\varphi , \psi]$ denote $[\varphi ,\psi] _0$.
  The graded $S$-bimodule
  $D^{n, G}_q $ or simply $D^n$ (where there is no confusion of the underlying
  group) is inductively defined as follows:
  We let $D^n = 0 $ for $n < 0$.  For each $n \geq 0$, we let $D^n$ be the
  $S$-bimodule generated by the set
  \begin{align*}
    Z_{n,q} = \{ \text{homogeneous } \varphi \in 
    \grHom_{\K}(S,S) \mid \text{ there is some $g \in G$ such that}\\
    [\varphi, r]_g \in D^{n-1} 
    \text{ for all homogeneous \(r \in S\)}  \}.
  \end{align*}
  The $R$-bimodule
  $D_q^G (S) = \bigcup D^n$ is a ring since $D^i D^j \subset D^{i+j}$.
  This implies that $D^0$ is a ring, and each $D^n$ is a $D^0$-module.

  For our ring $R$, we consider the bicharacter
  \begin{align*}
    \beta : \Gamma \times \Gamma &\to \K^*,\\
    ((a,b),(m,n)) &\mapsto q^{\frac{am + bn}{2}}.
  \end{align*}
   
  The \term{grading action} of $\Gamma$ on $R$ is the map 
  which assigns to $(a,b) \in \Gamma$ the automorphism of $R$
  \[
  \sigma _{(a,b)} (x^ny^m) = \beta((a,b),\deg(x^ny^m)) x^ny^m = 
  q^{\frac{a(n-m)+b(n+m)}{2}}x^ny^m.
  \]
  The set of all such grading homomorphisms will be denoted 
  $\sigma_\Gamma$.  These are always $q$-differential operators of 
  order 0.  In order to make this exposition more readable,  
  we will sometimes write $\sigma_x$ for $\sigma_\mba$ and 
  $\sigma_y$ for $\sigma_\mbb$.
  
  \begin{definition}
    For each $r \in R$, we define two multiplication homomorphisms 
    $\lambda_r$ and $\rho_r$ 
    in $\grhom$  by
    \begin{align*}
      \lambda_r(s) &= rs,\\
      \rho_r(s) &= sr.
    \end{align*}
  \end{definition}
  When a ring is commutative, these two homomorphisms coincide.
  The ring $R$ is very 
  close to being commutative in the following sense:
  \begin{align*}
    \lambda_x &=  \rho _x \sigma_y,\\
    \lambda_y &= \rho_y \sigma_x^{-1}.
  \end{align*}
  Since $\sigma_\Gamma$ and these multiplication maps generate all
  $q$-differential operators of order 0, we have the following lemma.
  \begin{lemma}
    The $q$-differential operators of order 0 on $\Gamma$-graded $R$,
    are generated by
    $\sigma_{\Gamma}$, $\lambda_x$, and $\rho _y$.
  \end{lemma}
  We have a corresponding Lemma for the $\Lambda$-graded ring:
  the ring $D_q^{0,\Lambda}$ is generated by 
  $\sigma _{\Lambda}$, $\lambda_x$, and $\rho_y$.  This leads us 
  to the general statement that 
  for each $n\geq 0$, we have 
  \[
     D^{n,\Lambda} \subset D^{n,\Gamma}.
  \]
  
  Before proceeding, we need some notation.  For any integers 
  $n$ and $a$ where $a \neq 0$, we put
  \[
     [ n ]_a = \frac{q^{na} - 1}{q^a-1} = 1 + q^a + q^{2a} + \cdots 
     + q^{(n-1)a}.
  \] 
  Denote $[n]_1$ simply by $[n]$, and define $[n]_0 = n$.  
\section{The first order $q$-differential operators}\label{S:D1}

  We will define a family of homomorphisms which fill the role of 
  the derivations on $R$ in the context of ordinary differential operators.
  \begin{definition}
    Let $a$ and $b$ be integers.  Define $\partial^{\beta^a}_x$
    and $\partial^{\beta^b}_y$ by 
    \begin{itemize}
    \item[]     $\partial ^{\beta^a}_x(x^iy^j) = [i]_a x^{i-1}y^j$,
    \item[]     $\partial ^{\beta^b}_y(x^iy^j) = [j]_b x^iy^{j-1}$.
    \end{itemize}
    When either $a$ or $b$ is 1, we will omit it as a superscript.
    We will simply write $\partial_x$ and $\partial_y$ for 
    $\partial^{\beta^0}_x$ and $\partial^{\beta^0}_y$ respectively.
  \end{definition}

  We note the following
  \begin{align*}
    [\partial ^{\beta^a}_x , \rho _y] &=0, 
       &[\partial ^{\beta^a}_x , \lambda _x] &=\sigma_x^a,\\
    [\partial ^{\beta^b}_y , \lambda _x] &=0, 
       &[\partial ^{\beta^b}_y , \rho _y] &=\sigma _y^b.
  \end{align*}
  Here, $\sigma _x ^a, \sigma _y^b$ denote $(\sigma _x )^a, (\sigma _y)^b$
  respectively.
  This shows that
  $\partial^{\beta^a}_x$ and $\partial^{\beta^b}_y$
  are $q$-differential operators of order 1.
  Furthermore, for $a,b \geq 2$
  \begin{align*}
    \partial _x^{\beta ^a} 
      &= \left( \frac{1-q}{1-q^a} \right) \partial _x^{\beta }
                (1 + \sigma_x +\cdots + \sigma_x^{a-1}),\\
    \partial _y^{\beta ^b} 
      &= \left( \frac{1-q}{1-q^b} \right) \partial _y^{\beta }
                (1 + \sigma_y +\cdots + \sigma_y^{b-1}).
  \end{align*}
  For $a,b \leq -1$, we have
  $\partial_x^{\beta^a} = \sigma_x^{a} \partial_x^{\beta^{-a}}$
  and
  $\partial_y^{\beta^b} = \sigma_y^{b} \partial_x^{\beta^{-b}}$.
  So with only $\K$, $\sigma_\Lambda$ and the four homomorphisms
  $\partial_x$, $\partial_y$, $\partial_x^{\beta}$, and 
  $\partial_y^{\beta}$, we may generate all other 
  $\partial_x^{\beta^{a}}$'s and $\partial_y^{\beta^{b}}$'s.
  We now prove the main theorem of this section.
  \begin{theorem}\label{T:D1}
    The $D^{0,\Gamma}$-module $D^{1,\Gamma}$,
    of first order $q$-differential operators 
    is generated by 
    $\{ 1$, $\partial^{\beta}_x$, $\partial^{\beta}_y$, $\partial_x$, 
    $\partial _y \}$.
  \end{theorem}
  \begin{proof}
    In this proof, we let $D^0 = D^{0, \Gamma}$ and $D^1 = D^{1, \Gamma}$.
    The $R$-bimodule $D^1$ of first order $q$-differential operators
    is generated as a $D^0$-module by homogeneous $\varphi \in \grhom$
    such that $[\varphi, \lambda_r] \in D^0$ for any $r\in R$ (see
    Corollary~1.2.1 of \cite{IM}).  Since $R$ is generated by the 
    two elements $x$ and $y$, it is enough to consider
    those homogeneous $\varphi$ such that $[\varphi ,\lambda_x]$ and
    $[\varphi ,\lambda_y]$ are in $D^0$.
    We need to show that $\varphi$ is
    in the $D^0$ span of the $\partial_x^{\beta^a}$'s and 
    $\partial_y^{\beta^b}$'s.
    
    Since $\rho_{\varphi(1)}$ is in $D^0$, we can replace $\varphi$ with 
    $\varphi-\rho_{\varphi(1)}$ and so assume that $\varphi(1) = 0$.
    We consider the following cases separately:\hfill
    \begin{enumerate}
\item If $\varphi$ is a 
      first order $q$-differential operator such that
      \begin{align*}
        [\varphi , \lambda_x] &= 0,\\
        [\varphi ,\lambda_y] &= \sigma _{(c,d)},
      \end{align*}
      then $\varphi (xy) = q \varphi (yx)$ implies that
      $d = -c-2$.  Call such pairs $(c,-c-2)\in \Gamma$ 
      type C1.
      Since $\varphi (1) = 0$,
      we can conclude 
      \[
        \varphi (x^n y^m) = [m]_b x^ny^{m-1},
      \]
      where $b = \frac{-c + d}{2} = -(c+1)$.
      Hence $\varphi = \partial _y ^{\beta^b}$.
      
      Similarly, if $\varphi$ is a 
      first order $q$-differential operator such that
      \begin{align*}
        [\varphi , \lambda_x] &= \sigma _{(c,d)},\\
        [\varphi ,\lambda_y] &=0,
      \end{align*}
      then $\varphi (xy) = q \varphi (yx)$ implies that
      $d = c+2$.  Call such pairs $(c,c+2) \in \Gamma$ 
      type C2.
      Here again
      we can conclude 
      \[
        \varphi (y^nx^m) = [m]_a y^n x^{m-1}
      \]
      where $a =c+1$. 
      Thus we have
      $\varphi = \partial^{\beta^a}_x$.

\item Let $\varphi$ be homogeneous such that  
      \begin{align*}
        [\varphi, \lambda_x] &= 0 ,\\
        [\varphi, \lambda_y] &= \rho _s \sum _{i} \alpha _i\sigma _{(c_i,d_i)},
      \end{align*}
      for $\alpha _i \in \K$.  As $[\rho_s, \lambda_r] =0$ for $r\in
      R$, we have $[\rho_s \varphi', \lambda_r]=\rho_s
      [\varphi',\lambda_r]$ for any $r \in R$ and $\varphi' \in \grhom$.  
      Hence, we can assume without loss of generality that
      $s=1$.  This implies that
      \[
        \varphi (x^ny^m) = \left( \sum _i \alpha _i
                              [m]_{b_i} \right) x^ny^{m-1}
      \]
      where $b_i = {\frac{-c_i +d_i}{2}}$.
      It follows that
      $[\varphi ,\lambda_y] (x^ny^m) = \frac{1}{q^n} (\sum _i \alpha _i q^{m
              b_i})x^ny^m$.
      Now we can use the fact that
      $[\varphi ,\lambda_y ] = \sum _{i} \alpha _i\sigma _{(c_i,d_i)}$
      to obtain
      \[
        \frac{1}{q^n} \sum _i \alpha _i q^{m
              b_i} =
              \sum _i \alpha _i q^{n(\frac{c_i +d_i}{2})} 
                      q^{mb_i}.
      \] 
      Let     
      $A_a = \{ i \mid \frac{c_i + d_i}{2} = a \textit{ and } \alpha _i 
      \neq 0\}.
      $
      We have,
      \[
        \sum_a \left(\sum_{i \in A_a} \alpha_i q^{mb_i}
                      \right)         
                      = \sum _a q^{n(a +1)}
                      \left( \sum _{i\in A_a} \alpha _i 
                              q^{mb_i}\right )
      \]
      for all nonnegative integers $m$ and $n$.
      This further implies that
      \[
        \sum _{a\neq -1} 
              \left(\sum _{i \in A_a} \alpha _i q^{mb_i}
                      \right)         
                      = \sum _{a\neq -1} q^{n(a +1)}
                      \left( \sum _{i \in A_a} \alpha _i 
                              q^{mb_i}\right )
      \]
      for all $m,n\in \mathbb{N}$.  If $A_a$ is not  
      empty for $a \neq -1$,
      then we can fix an $m$ such that
      $\gamma _a =\sum _{i\in A_a} \alpha _i q^{mb_i} 
      \neq 0$
      for some $a \neq -1$.
      Note that if $c_i + d_i = c_j +d_j$ and
      $(c_i,d_i) \neq (c_j ,d_j)$, then
      $-c_i +d_i \neq -c_j + d_j$.  This implies that either
      $a = -1$ or $\alpha _i = 0$ for
      $i \in A_a$ which contradicts our assumption.
      Hence, 
      \begin{align*}
        [\varphi ,\lambda_x] &= 0 ,\\
        [\varphi ,\lambda_y] &= \sum _{i} \alpha _i\sigma _{(c_i,d_i)},
      \end{align*}
      where $c_i +d_i = -2$ (that is, type C1).  
      Such a $\varphi$ is sum of first order
      differential operators covered in the first case.
      
      The case where
      $\varphi$ is a homogeneous homomorphism such that       
      \begin{align*}
        [\varphi ,\lambda_y] &= 0 ,\\
        [\varphi ,\lambda_x] &= \sum _{i} \alpha _i\sigma _{(a_i,b_i)},
      \end{align*}
      for $\alpha _i \in \K$, is dealt very similarly to the case just 
      discussed.
\item Let $\varphi$ be homogeneous such that
      \begin{align*}
        [\varphi,\lambda_x] 
               &= \rho _r \sum _{i} \alpha _i\sigma _{(a_i,b_i)} \\
        [\varphi ,\lambda_y] 
               &= \rho _s \sum _{j} \gamma _j\sigma _{(c_j,d_j)},
      \end{align*}
      for $\alpha _i,\gamma _j \in \K$.
      Since $\varphi$ is homogeneous, we have $r = x^{n+1}y^m$ and
      $s = x^n y^{m+1}$ for some integers $m$ and $n$.  
      This implies that $x$ is a factor of $r$.
      Whenever $a_i + b_i \neq 0$, we can rewrite
      $\rho_r \sigma _{(a_i, b_i)}$ as
      $\rho_{r^{\prime}} x \sigma _{(a_i +1,b_i-1)}$ for
      some $r^{\prime} \in R$ 
      (because $\lambda _x = \rho _x \sigma_y$).
      Since
      \[
        [\sigma _{(a,b)}, \lambda_x] = (q^{\frac{a+b}{2}} -1)x\sigma _{(a,b)},
      \]
      we have
      \[
        (q^{\frac{a_i+b_i}{2}}-1)^{-1}
        [\rho_{r'}\sigma _{(a_i+1,b_i-1)},\lambda_x]
                           = \rho_r\sigma _{(a_i,b_i)}.
      \]
      Hence, we can assume that $\varphi$ satisfies
      \begin{align*}
        [\varphi, \lambda_x] 
                &= \rho _r \sum _{i} \alpha _i\sigma _{(-a_i,a_i)} \\
        [\varphi, \lambda_y] 
                &= \rho _s \sum _{j} \gamma _j\sigma _{(c_j,d_j)},
      \end{align*}
      where $r = x^{n+1}y^m$ and
      $s = x^n y^{m+1}$.
      The relation
      \[
        [\varphi, \lambda_x] = \rho _r \sum _{i} \alpha _i\sigma _{(-a_i,a_i)}
      \]
      implies that
      \[
        [\varphi ,\lambda_x^n] =  nx^{n-1}\rho _r \sum _{i} \alpha _i
                        \sigma _{(-a_i,a_i)}.
      \] 
      Since $\varphi (1) =0$, 
      we obtain
      \[
        \varphi (x^t y) = \left( \sum _j \gamma _j + \frac{t} {q^{n+1}}
                                \sum _i \alpha _i q^{a_i} \right) 
                              x^{t+n}y^{m+1}
      \]                                 
      for all $t \geq 0$;
      whereas,
      \[
        \varphi (q^tyx^t ) = \left( \frac{t}{q^n}\sum _i \alpha _i + 
                                \sum _j \gamma _j 
                                q^{(\frac{c_j +d_j}{2} +1)t}
                                \right) x^{t+n}y^{m+1},
      \]
      for all $t\geq 0$.
      This implies that for all $t \geq 0$,
      \begin{align*}
        \sum _j \gamma _j &= \sum _j \gamma _j
                                           q^{(\frac{c_j +d_j}{2} +1)t}
                                             \text{ and },\\
        \frac{1}{q^n}{\sum _i \alpha _i} &= 
        \frac{1}{q^{n+1}}{\sum _i \alpha _i q^{a_i}}.   
      \end{align*}
      Hence the pairs $(c_j,d_j)$ are of type $C2$
      or  $\gamma_j =0$ for all $j$.  Either conclusion 
      takes us back to the second case.
    \end{enumerate}
    This proves the theorem.
  \end{proof}

\section{The generators of $\D$}\label{S:Dlambda}

  Let us define two subalgebras of $\grhom$, 
  \[
  D_x = \K \langle \lambda _x ,\; \partial_x ^{\beta^a} \mid a = -1,0,1 \rangle
  \text{\quad and \quad}
  D_y = \K \langle \rho _y ,\; \partial_y^{\beta ^b} \mid b = -1,0,1 \rangle .
  \]
  Since $[ \partial_x^{\beta^a}, \lambda_x]= \sigma_x^a$, 
  and $[\partial_y^{\beta^b}, \rho_y] = \sigma_y^b$, we have
  $D_x$ to be the $\K$-subalgebra generated by $\lambda_x$
  $\partial^{\beta}_x$, $\partial _x$, and $\sigma_x^{-1}$.
  A similar statement can be made for $D_y$.
  Each of $D_x$ and $D_y$ 
  is isomorphic to the ring of quantum differential operators on a 
  polynomial ring with one variable (see \cite{IM}).
  Moreover, elements of $D_x$ commute with those of $D_y$.

  The rings $D_x$ and $D_y$ are subrings of $\D$.
  We will show that $D_x$ and $D_y$ generate all of $\D$.
   
  The following is a straight-forward generalisation of Lemma~2.0.2 in 
  \cite{IM}.
  \begin{lemma}\label{lemma:0}
    For any $\varphi \in \D$ and any $r \in R$, there are 
    $\mb{c_1}$, $\mb{c_2}$, $\ldots$, $\mb{c_k}$ in $\Lambda$ such that
    \[
      [ \cdots [\ [ \varphi, \lambda_r]_\mb{c_1}, \lambda_r]_\mb{c_2} \cdots, 
      \lambda_r]_\mb{c_k} = 0.
    \]
  \end{lemma}
   
  \begin{corollary}
    For $r=x$, the $\mb{c_i}$ may be taken to be multiples of $\mba$.
  \end{corollary}
  \begin{proof}
    For any $\mb{c} \in \Lambda$ we have,
    $\beta(\mb{c}, (1,1)) = q^m$ for some integer $m$.  Thus,
    \[ 
      [\varphi, \lambda_x]_\mb{c} = 
      \varphi \lambda_x - q^m \lambda_x \varphi = 
      [\varphi, \lambda_x]_{(m,m)}
    \]
    for any graded homogenous endomorphism $\varphi$.
  \end{proof}
   
  \begin{lemma}\label{L:freeD0}
    The operators 
    $\{ \lambda _x^i, \rho _y^j, \sigma _{\mb{d}} \}_{i,j\geq 0,\mb{d}\in
      \Lambda } $ are 
    linearly independent over $\K$.
  \end{lemma}
  \begin{proof}
    Since the given collection of operators are graded, it is enough
    to check for linear independence of homogeneous operators.  This 
    reduces to check that the set $\{ \sigma _{\mb{d}} \}$ is linearly
    independent over $\K$ which is clear.
  \end{proof}

  \begin{lemma}\label{lemma:tims1}
    If $\varphi \in \D$ and $[\varphi, \rho_y]=0$, 
    then $\varphi \in D_x \langle \rho_y \rangle$.
  \end{lemma}
  \begin{proof}
    By the corollary above there are $\mb{c_1}$, $\mb{c_2}$, $\ldots$, 
    $\mb{c_k}$ which are multiples of $\mba$ and such that 
    \[
      [ \cdots [[ \varphi, \lambda_x]_\mb{c_1}, \lambda_x]_\mb{c_2} \cdots, 
      \lambda_x]_\mb{c_k} = 0.
    \]
    Put $\varphi_1 = \varphi$ and 
    $\varphi_{i+1} = [\varphi_{i},\lambda_x]_\mb{c_i}$.
    The operators $[\cdot, \lambda_x]_\mb{c_i}$ and $[\cdot, \rho_y]$
    commute, so we have $[\varphi_i, \rho_y]=0$ for each $i$.  We will 
    prove the lemma by descending induction on $k$.
    
    To start, $[\varphi_{k}, \lambda_x]_\mb{c_k} = 0$, so $\varphi_{k}$
    is in $D^{0,\Lambda }$. Hence 
    \[ \varphi_{k}= \sum_{i,j\in\Z,\mb{d}\in\Lambda} \alpha_{ij\mb{d}}
       \rho_{y}^i\lambda_{x}^j\sigma_\mb{d}
    \]
    where $\alpha_{ij\mb{d}} \in \K$.
    However, $[\varphi_{k}, \rho_y] = 0$, so for each $i$, $j$, and 
    $\mb{d}$, either $\alpha_{ij\mb{d}} = 0$ or $\sigma_\mb{d}$ 
    commutes with $\rho_y$ (by Lemma~\ref{L:freeD0}).  
    Hence, each $\alpha_{ij\mb{d}}\sigma_\mb{d}$ 
    is in $D_x$, and thus $\varphi_{k}$ is in $D_x\langle \rho_y\rangle$.
    
    Now suppose that $\varphi_{i}$ is in $D_x\langle \rho_y \rangle$.
    Then so is $\sigma^{-1}_{\mb{c_{i-1}}} \varphi_{i}$ 
    by our assumptions on the 
    $\mb{c_i}$.
    Then we have
    \[
      \sigma^{-1}_{\mb{c_{i-1}}}\varphi_{i}
      = [\sigma^{-1}_{\mb{c_{i-1}}}\varphi_{i-1}, \lambda_x] = 
      \sum_j f_j \rho_{y}^j
    \]
    for some $f_j \in D_x$.
    By \cite{IM}, we have $F_j$ in $D_x$ such that 
    \[ [F_j, \lambda_x] = f_j. \]
    Put 
    \[ \psi = \sum F_j \rho_{y}^j. \]
    Then $\psi \in D_x\langle\rho_y\rangle$, and we have 
    \begin{align*}
       [\psi - \sigma^{-1}_{\mb{c_{i-1}}}\varphi_{i-1}, \lambda_x] 
       = & [\psi, \lambda_x] 
             - [\sigma^{-1}_{\mb{c_{i-1}}}\varphi_{i-1}, \lambda_x ] \\
       = & \sum f_j \rho_{y}^j
             - \sigma^{-1}_{\mb{c_{i-1}}}
             [\varphi_{i-1}, \lambda_x]_\mb{c_{i-1}}  \\
       =& 0 .
    \end{align*}
    As $\sigma^{-1}_{\mb{c_{i-1}}}$, $\varphi_{i-1}$, and $\psi$ all 
    commute with $\rho_y$, 
    $[\psi - \sigma^{-1}_{\mb{c_{i-1}}}\varphi_{i-1}, \rho_y] = 0$.
    This implies that $\psi - \sigma _{\mb{c_{i-1}}}^{-1} \varphi _{i-1}
    \in D^{0,\Lambda}$.
    By applying the base case to 
    $\psi - \sigma^{-1}_{\mb{c_{i-1}}} \varphi_{i-1}$, 
    we conclude  $\psi - \sigma^{-1}_{\mb{c_{i-1}}} \varphi_{i-1} 
    \in D_x \langle \rho _{y} \rangle$.
  \end{proof}
  \begin{remark}
    The lemma above is interesting in its own right.  A similar proof can be
    used to show that if $\varphi \in \D$ and $[\varphi, \lambda _x]=0$, 
    then $\varphi \in D_y \langle \lambda _x \rangle$.
  \end{remark}
  \begin{theorem}
    The ring $\D$ is generated by $D_x$ and $D_y$.
  \end{theorem}
  \begin{proof}
    Let $\varphi$ be a homogeneous $q$-differential operator of order 
    $m$.  Then $\varphi = \sum \varphi_i$ where each $\varphi_i$ has the 
    property that for some $\mb{c}\in \Lambda$, 
    $[\varphi_i, \lambda_r]_\mb{c}$ is in $D^{m-1,\Lambda}$ for any
    $r \in R$.  Then $\varphi$ will be in the algebra generated by 
    $\partial_x$, $\partial_y$, $\partial_x^{\beta}$, 
    $\partial_y^{\beta}$, and $D^{0,\Lambda}$ if each $\varphi_i$ is so.
    Thus we will assume that there is a $\mb{c} \in \Lambda$ such
    that for any $r \in R$, $[\varphi, \lambda_r]_\mb{c} \in D^{m-1,\Lambda}$
    and proceed by induction on $m$.
    
    The base case is immediate, so we will suppose already that 
    all $q$-differential operators of order $m-1$ or less 
    can be expressed with the generators in the set above.
    In particular, $[\varphi,\lambda_y]_\mb{c}$ is in the span of 
    $D_x$ and $d_y$.
    Since $[\varphi,\lambda_y]_\mb{c} = [\varphi,\rho_y]_\mb{d} \sigma_x$
    for some $\mb{d} \in \Lambda$ which depends on $\mb{c}$ and $\varphi$, 
    we have $[\varphi, \rho_y]_\mb{d}$ is in the span of $D_x$ and $D_y$.
    Since $[\sigma_\mb{d}^{-1}\varphi,\rho_y] 
    = \sigma_\mb{d}^{-1}[\varphi, \rho_y]_\mb{c}$, 
    and $\mb{c}\in \Lambda$, we can assume without loss of 
    generality that $[\varphi, \rho_y] = \sum f_i g_i$ for some 
    $f_i \in D_x$ and $g_i \in D_y$.
    
    By Theorem~2.0.1 of \cite{IM}, 
    there are $G_i \in D_y$ such that $[G_i, \rho_y] = g_i$.
    Put $\psi = \sum f_i G_i - \varphi$.
    Then $[\psi, \rho_y] = 0$.  Hence $\psi \in D_x\langle\rho_y\rangle$. 
    It follows that $\varphi$ is in the ring generated by $D_x$ and $D_y$.
  \end{proof}

\section{Properties of the ring $\D$}\label{S:pro}
  In this section, we consider the properties of $\D$ as a ring.
  \begin{proposition}\label{P:tensor}
    The ring $D_x \bigotimes _{\K} D_y$ is simple.
  \end{proposition}
  \begin{proof}
    Let $I$ be an ideal of $D_x \bigotimes_{\K} D_y$.
    Every element $\theta$ of $I$ can be written as
    $\sum_i \psi_i \otimes \eta _i$ where $\psi _i \in D_x$ and 
    $\eta _i \in D_y$ for all $i$. Since 
    $\theta (\lambda_x\tensor 1) - q^n (\lambda_x \tensor 1) \theta$
    is in $I$ for every integer $n$, we conclude there is some 
    $\theta \in I$ such that each $\psi_i$ is in $\K[\lambda_x]$.
    Similarly, we can assume that each $\eta_i$ is in $\K[\rho_y]$.
    When we consider the commutators of $\theta$ with 
    $\partial_x \otimes 1$ and $1 \otimes \partial_y$, we conclude 
    $I$ contains $1\otimes 1$.
  \end{proof}
  \begin{theorem}
    The ring $\D$ is isomorphic to $D_x \bigotimes_{\K} D_y$ as graded 
    rings.
  \end{theorem}
  \begin{proof}
    Note that the elements of $D_x$ and those of $D_y$ commute with
    each other and
    that the algebras $D_x$ and $D_y$ generate $\D$.
    Hence, we have a surjective map from 
    $D_x \bigotimes _{\K} D_y$ to $\D$.  This map is injective by the 
    proposition above.  Hence the theorem.
  \end{proof}
  We fix the following notations: \hfill
  \begin{enumerate}
    \item[]       $\tau _x = \lambda _x \partial _x$.
    \item[]       $\tau _y = \rho _y \partial _y$.
  \end{enumerate}
  \begin{theorem}
    The ring $\D$ is a domain.
  \end{theorem}
  \begin{proof}
    For each $\mb{c} \in \Lambda$, let
    $(\D)_{\mb{c}}$ denote the set of all homogeneous operators in 
    $\D$  of degree $\mb{c}$. Then 
    $(\D)_{(0,0)} = (D_x)_0 \bigotimes _{\K} (D_y)_0$
    where $(D_x)_0$ and $(D_y)_0$ are respectively the ring of 
    operators in $D_x$ and $D_y$ of degree 0.  
    By Lemma~3.0.5 of \cite{IM}, 
    $
    (D_x)_0 = \K [ \tau _x, \sigma _x ^{\pm 1}]$ and $
    (D_y)_0 = \K [\tau _y, \sigma _y ^{\pm 1}]
    $
    are localised polynomial rings.
    Hence, 
    $
    (\D)_{(0,0)} = 
    \K[ \tau _x, \tau_y, \sigma _x^{\pm}, \sigma _y^{\pm}] 
    $
    is a localised polynomial ring and in particular
    a domain. 
    
    The operator $\lambda _x$ is not a left zero-divisor in 
    $\D$ because the ring $R$ is a domain.
    Similarly, $\rho_y$ is not a right zero-divisor.
    Suppose that $\lambda_x$ is a right zero-divisor in $\D$ and 
    $\varphi \lambda_x = 0$.  Then for any $\mb{c} \in \Lambda$, 
    $[\varphi, \lambda_x]_{\mb{c}}= q^m \lambda_x \varphi$
    where $\sigma_{\mb{c}}(x) = q^m x$.  It follows that for 
    any $\mb{c_1}$, $\mb{c_2}$, $\ldots$, $\mb{c_k}$ in $\Lambda$ 
    there is some integer $m$ such that
    \[
       [\cdots [\ [ \varphi, \lambda_x]_\mb{c_1}, \lambda_x]_\mb{c_2} 
        \cdots, \lambda_x]_\mb{c_k} = q^m \lambda_x^k \varphi.
    \]
    Since $\lambda_x$ is not a left zero-divisor, this cannot be
    0.  Hence, by the Lemma~\ref{lemma:0}, $\varphi$ is not in 
    $\D$.  A similar argument shows that $\rho_y$ cannot be a 
    zero-divisor on the left.
    
    Finally, suppose $\varphi \psi = 0$.  It is sufficient to consider 
    the case when $\varphi, \psi$ are homogeneous of degrees 
    $(a,b),(c,d)$ respectively, in the sense that 
    $\varphi (x^n y^m) = \alpha x^{n+a}y^{m+b}$,
    and similarly for $\psi$.
    By suitably multiplying by (or factoring of) powers of 
    $\lambda _x$ and $\rho_y$,
    we can assume that $\varphi$ and $\psi$ are in $(\D)_{(0,0)}$, 
    which we know is a domain. Hence the theorem.
  \end{proof}


\section{The motivating diagram}\label{S:mot}
  The universal enveloping algebra $U(\slin)$ of the Lie algebra 
  $\slin$ acts on its representation ring $\K[u,v]$ through 
  $D(\K[u,v])$, the ring of differential operators on $\K[u,v]$.  
  Similarly, the quantised universal enveloping algebra (or 
  \textit{quantum group}) $U_q(\slin)$, a deformation of $U(\slin)$, 
  acts through the quantum differential operators of its 
  representation ring, $R$.
  Understanding these actions on the underlying representation rings 
  is the goal of this section.

\subsection{Quantum group on $\mathbf{\slin}$}
  Let $U_q$ denote the quantum group corresponding to the Lie algebra
  $sl_2(\K)$.  
  That is, $U_q$ is a $\K$-algebra generated by $E,F,K,K^{-1}$,
  with relations given by
  \begin{align*}
    KK^{-1} = &1 = K^{-1}K, \\
    KEK^{-1} &= q^2E,\\
    KFK^{-1} &= q^{-2}F,\\
    EF - FE &= \frac{K - K^{-1}}{q-q^{-1}}.
  \end{align*}
  Extensive treatments of $U_q$ can be found in \cite{CP} and \cite{J},
  but we only need from these the fact that $U_q$ is a Hopf-algebra
  with a co-multiplication map $\Delta$ given by
  \begin{align*}
    \Delta (E) &= E\otimes 1 + K \otimes E,\\
    \Delta (F) &= F \otimes K^{-1} + 1 \otimes F,\\
    \Delta (K) &= K \otimes K.
  \end{align*} 

\subsection{The integral form of $\mathbf{U_q}$}
  Let $\CA = \Q[q, q^{-1}]$ be the Laurent polynomial ring over $\Q$.   
  Let $\dbox{m}= \frac{q^m - q^{-m}}{q -q^{-1}}$
  and $\dbox{m} ! = \prod_{1\leq i \leq m}\dbox{i}$. 
  By convention, $\dbox{0} ! =1$.
  The integral form $U_{q,\CA}$ of 
  $U_q$ is the $\CA$-subalgebra
  of $U_q$ generated by $E^{(m)},F^{(m)}, K , K ^{-1}$ for
  positive integers $m$ where
  \[
      E^{(m)} =\frac{E^m}{\dbox{m}!},\qquad F^{(m)} = \frac{F^m}{\dbox{m}!}.
  \]

\subsection{The ring $\mathbf{D_{\CA}}$ }
      Let $D_{\CA}$ be the $\CA$-subalgebra of
      $D_q^{\Lambda}(R)$ generated by
      \[
      \{ \lambda _x, \rho _y , \sigma _{\mb{d}} , \partial _x, \partial _y,
      (\partial _x^{\beta})^{(m)}, (\partial _y^{\beta})^{(m)} \}_{m\in 
        \mathbb{N},
        \mb{d}\in \Lambda}
      \]
      where
      \[
      (\partial _x^{\beta })^{(m)} = 
                      \frac{(\partial _x^{\beta})^m}{[2]^m\dbox{m}!},
                              \qquad
      (\partial _y^{\beta})^{(m)} =
                      \frac{(\partial _y^{\beta})^m}{[2]^m\dbox{m}!}.
      \]
      The completion of localisation of $D_{\CA}$ at $(q-1)$ 
      will be denoted by $\widehat{D_\CA}$.

\subsection{The action of $\mathbf{U(sl_2)}$ on $\mathbf{\Q[u,v]}$.}
      The ring of $\Q$-linear usual differential operators on
      the polynomial ring $\Q[u,v]$ is the second Weyl algebra.
      That is, $D(\Q[u,v])$ is a $\Q$ algebra 
      with generators $\lambda _u ,\lambda _v, \partial _u, \partial _v$
      and relations $[\partial _u, \lambda_u] = 1,
      [\partial _v, \lambda_v] = 1$, and all the other commutators equal 0.
      
      The action of $U(\slin)$ on $\Q[u,v]$ gives a 
      homomorphism $\psi : U(\slin) \to D(\Q[u,v])$ by
      \begin{align*}
              \psi(e) &= \lambda _u\partial_v,\\
              \psi(h) &= \lambda _u\partial_u - \lambda _v\partial_v,\\
              \psi(f) &= \lambda _v\partial_u.
      \end{align*}
      Let $\CA [ \bar{h}]$ denote a polynomial ring over $\CA$ in one
      variable $\bar{h}$.  We let
      \[
      \psi (\bar{h} ) = \lambda _u\partial_u + \lambda _v\partial_v.
      \]
      Consider the projection of the centre  $Z \subset U(\slin)$ to 
      $\mathbb{Q} [ h ] = U^0 \subset U(\slin)$ corresponding to the
      triangular decomposition $U(\slin) = U^-U^0U^+$.
      Since $\mathbb{Q} [h] \cong \mathbb{Q} [\bar{h} ] $ as rings, 
      we can consider
      $\mathbb{Q} [\bar{h} ] $ as a $Z$-module via the projection.
      Then, the map $\psi$ gives rise to a map
      \[
      \psi:   U(\slin) \bigotimes _Z \mathbb{Q} [\bar{h}] \to 
                                D(\mathbb{Q} [u,v]).
      \]

      This map extends trivially to a homomorphism of 
      power series rings
      \[
       \psi : (U(\slin) \bigotimes _Z \mathbb{Q} [\bar{h}])[\![t]\!] 
                      \to D(\Q[u,v])[\![t]\!]
      \]
      which is $\Q[\![t]\!]$ linear.

\subsection{The action of $\mathbf{U_{q,A}}$ on $\mathbf{R}$.}
  The quantum group $U_q$ acts on the quantum plane as
  \begin{align*}
    K(1) &= 1,& K(x) &= qx,& K(y) &= \frac{1}{q}y,\\
    E(1) &= 0,& E(x) &= 0 ,& E(y) &= x,\\
    F(1) &= 0,& F(x) &= y ,& F(y) &=0
  \end{align*}
  with this action extended to all of $R$ via $\Delta$.  
  The resulting action can be listed as follows:
  \begin{align*}
    K(x^iy^j) &= q^{i-j}x^iy^j,\\
    E(x^iy^j) &= q^i[j]_{-2} x^{i+1}y^{j-1},\\
    F(x^iy^j) &= q^j[i]_{-2} x^{i-1}y^{j+1}.
  \end{align*}    
  That is, $K$ acts on $R$ as an automorphism, $E$ acts as a 
  left $K$-derivation, and $F$ acts as a right $K^{-1}$-derivation.

  The action of $E$ on $R$ is the same as the 
  $q$-differential operator 
  $\lambda_x \sigma_x \partial _y ^{\beta ^{-2}}$.
  This is an element of $D_\CA$ which can be expressed as  
  $\lambda_x \sigma_x \sigma_y^{-2} ({\partial_y^{\beta}})^{(1)}
  (1 + \sigma _y)$.  Similarly, we can identify 
  the action of $F$ on $R$ and with the action of an element of 
  $D_\CA$.  This leads us to define an $\CA$-linear homomorphism
  $
    \phi : U_{q,\CA } \to D_{\CA }
  $
  by putting
  \begin{align*}
    \phi(E) &= \sigma _{(2,0)} \lambda _x  (\partial_y^{\beta})^{(1)}
                   (1+\sigma_{(1,-1)}),\\ 
    \phi(F) &= \sigma _{(-2,0)} \rho _y (\partial_x^{\beta})^{(1)}
                   (1+\sigma_{(-1,-1)}),\\
    \phi(K) &= \sigma _{(2,0)}.
  \end{align*}

  If we define $d_y$ and $d_x \in D_{\CA}$ as
  \begin{align*}
    d_y &= \sigma _{(2,0)} (\partial_y ^{\beta})^{(1)}
                   (1 + \sigma_{(1,-1)}),\\ 
    d_x &= \sigma _{(-2,0)} (\partial_x ^{\beta})^{(1)}
                   (1 + \sigma_{(-1,-1)}),
  \end{align*}
  then  $d_y^{(m)} = d_y^m/[\![m]\!] !$ and $d_x^{(m)} = d_x^m/[\![m]\!]!$
  are elements of $D_\CA$.  Since $d_y \lambda _x = q \lambda _x d_y$
  and $d_x\rho_y = q \rho _y d_x$,
  we have
  \begin{align*}
    \phi (E^{(m)}) &= q^{\frac{m(m-1)}{2}}(\lambda_x)^m d_y^{(m)},\\
    \phi (F^{(m)}) &= q^{\frac{m(m-1)}{2}}(\rho _y)^m d_x^{(m)}.
  \end{align*}

  The action of $\sigma _{(0,2)}$
  is compatible with the action of the centre 
  $Z_q \subset U_{q,\CA }$.
  That is, consider the projection $p : Z_q \to U^o$
  corresponding to the triangular decomposition $U_q = U_q^-U_q^0U_q^+$.
  Now, think of the action of $K,K^{-1}$ as $\sigma _{(0,2)}$ and
  $\sigma _{(0,-2)}$ respectively.
  We can thus consider $\CA [\sigma _{(0,2)}]$ as a module over $Z_q$.
  Then the ring homomorphism 
  $\phi : U_{q,\CA} \bigotimes _{\CA} \CA [ \sigma_{(0,2)} ] \to D_{\CA }$
  factors through $Z_q$ to give a ring homomorphism
  \[
    \phi : U_{q,\CA} \bigotimes _{Z_q} \CA [ \sigma_{(0,2)} ] \to D_{\CA }.
  \]
  Let $\widehat{U_{q,\CA} \bigotimes _{Z_q} \CA [\sigma_{(0,2)}]}$ 
  denote the
  inverse limit of $U_{q,\CA} \bigotimes _{Z_q} \CA [ \sigma_{(0,2)} ]$ 
  with respect
  to the ideal $(q-1)$.
  Then we have the required map (by abuse of notation, we call it
  $\phi$ again)
  \[
     \phi : \widehat{U_{q,\CA} \bigotimes _{Z_q} \CA [ \sigma_{(0,2)} ]} \to 
                      \widehat{D_{\CA}}.
  \]

\subsection{The set-up}
      There is a commutative diagram 
      \[
      \begin{CD}
         \widehat{U_{q,\CA} \bigotimes _{Z_q}\CA [\sigma_{(0,2)} ]} @>{\phi}>> 
                              \widehat{D_\CA}\\
          @V{\mu}VV               @VV{\nu}V\\
         (U(\slin)\bigotimes _Z \CA [\bar{h}])[\![t]\!]  @>{\psi}>> 
                              D(\Q[u,v])[\![t]\!]
      \end{CD}
      \]
      which we shall now describe.

\subsection{The map $\mathbf{\nu}$.}
  For any expression $a$ in $D(\Q[u,v])$,  
  let 
  \[
      q^{(a)} = \sum_{n=0}^\infty \binom{a}{n} t^n,          
  \]
  where
  \begin{align*}
    \binom{a}{n} &= \frac{(a)(a-1)\cdots (a-n+1)}{n!}
         \textit{  for } n\in \mathbb{Z}, n\geq 1,\\
    \binom{a}{0} &= 1.
  \end{align*}
  We use parentheses to distinguish these power series in 
  $t$ from mere powers of $q$.  
  Then $q^{(1)} = 1+t$ and $q^{(a+b)}= q^{(a)}q^{(b)}$.

  For any $a \in D(\Q[u,v])$, the expression $q^{(a)} - 1$ 
  is divisible by $a$ and $t$.  Hence  
  \[
      P(a) = (q^{(a)} - 1) / a t
  \]
  is a well-defined power series in $t$ with coefficients 
  in $D(\Q[u,v])$, and it is invertible since its constant 
  coefficient is 1.
  
  We define the $\Q$-linear homomorphism
  $\nu: \widehat{D_\CA} \to D(\Q[u,v])[\![t]\!]$
  by  
  \begin{align*}
         \nu(q)           &= q^{(1)},&&\\
         \nu(\sigma_{x})  &= q^{(\lambda_u\partial_u)}
        &\nu(\sigma_{y})  &= q^{(\lambda_v\partial_v)},\\
         \nu(\lambda _x)  &= \lambda _u,
        &\nu(\rho_y)      &= \lambda _v,\\
         \nu(\partial _x) &= \partial _u,
        &\nu(\partial _y) &= \partial_v,\\
         \nu(\partial_x^{\beta}) &= \partial _u P(\lambda_u \partial _u),      
        &\nu(\partial_y^{\beta}) &= \partial _v P(\lambda_v \partial _v).
  \end{align*}
      
\subsection{The map $\mathbf{\mu}$.}
  Here we shall use the same notation as in the previous subsection.
  For $a \in \widehat{U_{q,\CA} \bigotimes _{Z_q} 
           \CA[ \sigma_{(0,2)}]}$,
  let $q^{(a)} = \sum _{n=0}^{\infty} \binom{a}{n} t^n$ 
  where the binomial coefficients are defined formally.
  Again, define the invertible $P(a) = (q^{(a)} - 1)/ a t$
  for any $a \in \widehat{U_q \bigotimes _{Z_q} \CA [\sigma_{(0,2)} ]}$. 
  Now we can define
  \begin{align*}
           \mu : \widehat{U_q \bigotimes _{Z_q} \CA [\sigma_{(0,2)} ]} 
                         &\to (U(\slin)\bigotimes _Z \mathbb{Q} 
                         [\bar{h}])[\![t]\!],\\
            \mu (q) &\mapsto (1+t),\\
            \mu (K) &\mapsto q^{(h)},\\
            \mu (E) &\mapsto \frac{ q^{(h)}e
                     P(\frac{-h + \bar{h}}{2}) 
                     (1 + q^{(\frac{h - \bar{h}}{2})})}{1+q^{(1)}},\\
            \mu (F) &\mapsto \frac{q^{(-h)} f
                     P(\frac{h + \bar{h}}{2})  
                     (1 + q^{(\frac{-h -\bar{h}}{2})})}{1+q^{(1)}},\\
            \mu (\sigma _{(0,2n)}) &\mapsto q^{(n\bar{h})}.
  \end{align*}
  This completes the commutative diagram.

\subsection{An Application}
  \begin{proposition}\label{C:UtoU_q}
    Given an irreducible, finite dimensional $U(\slin)$-module, 
    there is a natural way to define a $U_q$-action to give a 
    type-1 irreducible, finite dimensional module of $U_q$.
  \end{proposition}
  \begin{proof}
    Let $V$ be an $n+1$-dimensional irreducible module of $U(\slin)$, for
    $n\geq 0$.  
    Without loss of generality, we can assume
    that a $\K$-basis of $V$ is 
    \[
    \{ y^n, y^{n-1}x, y^{n-2}x^2,\cdots ,x^n \}
    \]
    where
    $x^n$ is the highest weight vector, and that the action of $U(\slin)$ 
    is given by the map $\psi$.   By using the map $\mu$, we get the 
    required $U_q$ action on $V$.
  \end{proof}


\section{Conclusions and Acknowledgements}
We suspect that $D^{\Gamma}_q (R) = D^{\Lambda}_q(R)\langle \sigma _{(1,0)}
\rangle$.  The evidence for this is that any left 
$\sigma_{(1,0)}$-derivation is a $q$-differential operator of order 0, 
not 1.  This suggests that if $\varphi$ is a $q$-differential operator
of order $m$, then so is $\varphi \cdot r - \sigma_{(1,0)}(r) \cdot \varphi$
for every $r$ in $R$.  It follows that $(D_q^{\Gamma}(R))^m$ would
be generated as a $(D_q^{\Gamma}(R))^0$-module by $Z_q^m$.  
 
Furthermore, we expect results similar to those in Section~\ref{S:mot}.
In particular, for a general semi-simple Lie algebra $\CG$ we 
expect a map analogous to $\mu$ could be constructed.

We thank Professor Valery~A.~Lunts for suggesting this question.
We also thank Northeastern Hill University in Shillong, Meghalaya, 
India for providing support and a beautiful working environment 
during the preparation of this paper.


\bibliographystyle{amsplain}
  
\end{document}